\input amstex
\documentstyle{amsppt}
%
\catcode`@=11
\redefine\output@{%
  \def\break{\penalty-\@M}\let\par\endgraf
  \ifodd\pageno\global\hoffset=105pt\else\global\hoffset=8pt\fi  
  \shipout\vbox{%
    \ifplain@
      \let\makeheadline\relax \let\makefootline\relax
    \else
      \iffirstpage@ \global\firstpage@false
        \let\rightheadline\frheadline
        \let\leftheadline\flheadline
      \else
        \ifrunheads@ 
        \else \let\makeheadline\relax
        \fi
      \fi
    \fi
    \makeheadline \pagebody \makefootline}%
  \advancepageno \ifnum\outputpenalty>-\@MM\else\dosupereject\fi
}
\def\Beta{\mathchar"0\hexnumber@\rmfam 42}
\catcode`\@=\active
\nopagenumbers
\chardef\textvolna='176
\def\negskp{\hskip -2pt}

\chardef\degree="5E
\def\blue#1{#1}

\catcode`#=11\def\diez{#}\catcode`#=6
\catcode`&=11\catcode`&=4
\catcode`_=11\def\podcherkivanie{_}\catcode`_=8
\catcode`~=11\def\volna{~}\catcode`~=\active
\def\mycite#1{\cite{\blue{#1}}\immediate\special{ps:
     ShrHPSdict begin /ShrBORDERthickness 0 def}}
\def\myciterange#1#2#3#4{\cite{\blue{#2#3#4}}\immediate\special{ps:
     ShrHPSdict begin /ShrBORDERthickness 0 def}}
\def\mytag#1{%
    \tag#1}
\def\mythetag#1{\thetag{\blue{#1}}\immediate\special{ps:
     ShrHPSdict begin /ShrBORDERthickness 0 def}}
\def\myrefno#1{\no#1}
\def\myhref#1#2{\blue{#2}\immediate\special{ps:
     ShrHPSdict begin /ShrBORDERthickness 0 def}}
\def\myEarXivlink{\myhref{http://arXiv.org}{http:/\negskp/arXiv.org}}

\def\mytheorem#1{\csname proclaim\endcsname{Theorem #1}}
\def\mytheoremwithtitle#1#2{\csname proclaim\endcsname{Theorem #1#2}}

\def\mylemma#1{\csname proclaim\endcsname{Lemma #1}}
\def\mylemmawithtitle#1#2{\csname proclaim\endcsname{Lemma #1#2}}

\def\mycorollary#1{\csname proclaim\endcsname{Corollary #1}}

\def\myconjecture#1{\csname proclaim\endcsname{Conjecture #1}}
\def\myconjecturewithtitle#1#2{\csname proclaim\endcsname{Conjecture #1#2}}

\def\myproblem#1{\csname proclaim\endcsname{Problem #1}}
\def\myproblemwithtitle#1#2{\csname proclaim\endcsname{Problem #1#2}}
\def\mytheproblem#1{\blue{#1}\immediate\special{ps:
     ShrHPSdict begin /ShrBORDERthickness 0 def}}

\font\eightcyr=wncyr8
\pagewidth{360pt}
\pageheight{606pt}
\topmatter
\title
On a pair of cubic equations associated with perfect cuboids.
\endtitle
\rightheadtext{On a pair of cubic equations associated with perfect cuboids.}
\author
Ruslan Sharipov
\endauthor
\address Bashkir State University, 32 Zaki Validi street, 450074 Ufa, Russia
\endaddress
\email\myhref{mailto:r-sharipov\@mail.ru}{r-sharipov\@mail.ru}
\endemail
\abstract
    A perfect cuboid is a rectangular parallelepiped with integer edges 
and integer face diagonals whose space diagonal is also integer. The existence 
of such cuboids is neither proved, nor disproved. A rational perfect cuboid
is a natural companion of a perfect cuboid absolutely equivalent to the latter
one. Its edges and face diagonals are rational numbers, while its space diagonal 
is equal to unity. Recently, based on a symmetry reduction, it was shown that 
edges of a rational perfect cuboid are roots of a certain cubic equation
with rational coefficients depending on two rational parameters. Face diagonals 
of this cuboid are roots of another cubic equation whose coefficients are
rational numbers depending on the same two rational parameters. In the present 
paper these two cubic equations are studied for reducibility. Six special cases 
of their reducibility over the field of rational numbers are found. 
\endabstract
\subjclassyear{2000}
\subjclass 11D41, 11D72, 12E05, 13A50\endsubjclass
\endtopmatter
\TagsOnRight
\document

%
%
\head
1. Introduction.
\endhead
     The problem of a perfect cuboid is known since 1719, but is still not resolved.
For the history of this problem the reader is referred to \myciterange{1}{1}{--}{44}.
Let $x_1$, $x_2$, $x_3$ be edges of a cuboid, $d_1$, $d_2$, $d_3$ be its face 
diagonals, and $L$ be its space diagonal. Then the cuboid is described by the following 
four polynomial equations:
$$
\xalignat 2
&\hskip -2em
x_1^2+x_2^2+x_3^2-L^2=0,
&&x_2^2+x_3^2-d_1^{\kern 1pt 2}=0,\\
\vspace{-1.7ex}
\mytag{1.1}\\
\vspace{1ex}
\vspace{-1.7ex}
&\hskip -2em
x_3^2+x_1^2-d_2^{\kern 1pt 2}=0,
&&x_1^2+x_2^2-d_3^{\kern 1pt 2}=0.
\endxalignat
$$
This paper continues the series of papers \myciterange{45}{45}{--}{50} applying the 
symmetry approach to the equations \mythetag{1.1}. Indeed, using three numbers
$x_1$, $x_2$, $x_3$, one can build the cubic equation $(x-x_1)(x-x_2)(x-x_3)=0$
which expands to
$$
\hskip -2em
x^3-E_{10}\,x^2+E_{20}\,x-E_{30}=0.
\mytag{1.2}
$$
Similarly, the equation $(d-d_1)(d-d_2)(d-d_3)=0$ expands to
$$
\hskip -2em
d^{\kern 1pt 3}-E_{01}\,d^{\kern 1pt 2}+E_{02}\,d-E_{03}=0.
\mytag{1.3}
$$
The coefficients $E_{10}$, $E_{20}$, and $E_{30}$ of the equation
\mythetag{1.2} are elementary symmetric polynomials of three variables 
$x_1$, $x_2$, $x_3$ (see \mycite{51}). Similarly, the coefficients $E_{01}$, $E_{02}$, 
and $E_{03}$ of the equation \mythetag{1.3} \pagebreak are elementary symmetric 
polynomials of three variables $d_1$, $d_2$, $d_3$. Here are the formulas for these
polynomials:
$$
\xalignat 2
&\hskip -2em
x_1+x_2+x_3=E_{10},
&&d_1+d_2+d_3=E_{01},\\
&\hskip -2em
x_1\,x_2+x_2\,x_3+x_3\,x_1=E_{20},
&&d_1\,d_2+d_2\,d_3+d_3\,d_1=E_{02},
\quad
\mytag{1.4}\\
&\hskip -2em
x_1\,x_2\,x_3=E_{30},
&&d_1\,d_2\,d_3=E_{03}.
\endxalignat
$$
Mixing $x_1$, $x_2$, $x_3$ with $d_1$, $d_2$, $d_3$, one can write the following
formulas:
$$
\hskip -2em
\aligned
&x_1\,x_2\,d_3+x_2\,x_3\,d_1+x_3\,x_1\,d_2=E_{21},\\
&x_1\,d_2+d_1\,x_2+x_2\,d_3+d_2\,x_3+x_3\,d_1+d_3\,x_1=E_{11},\\
&x_1\,d_2\,d_3+x_2\,d_3\,d_1+x_3\,d_1\,d_2=E_{12}.
\endaligned
\mytag{1.5}
$$
The left hand sides of the formulas \mythetag{1.4} complemented with the left hand 
sides of the formulas \mythetag{1.5} constitute the complete set of so-called
elementary multisymmetric polynomials. For the theory of multisymmetric polynomials,
either elementary and non-elementary, the reader is referred to 
\myciterange{52}{52}{--}{72}.\par
     The cuboid equations \mythetag{1.1} imply some equations for $E_{10}$, 
$E_{20}$, $E_{30}$, $E_{01}$, $E_{02}$, $E_{03}$, $E_{21}$, $E_{11}$, and $E_{12}$ 
in \mythetag{1.4} and \mythetag{1.5}. These equations are called factor equations.
They were studied in \mycite{46} and \mycite{47} using ideals in polynomial rings
and their Gr\"obner bases (see the general theory in \mycite{73}). In \mycite{48}
the factor equations were reduced to a single biquadratic equation for three
variables $E_{10}$, $E_{01}$, and $E_{11}$:
$$
\hskip -2em
(2\,E_{11})^2+(E_{01}^2+L^2-E_{10}^2)^2-8\,E_{01}^2\,L^2=0. 
\mytag{1.6}
$$
The other variables $E_{20}$, $E_{30}$, $E_{02}$, $E_{03}$, $E_{21}$, and $E_{12}$ 
are expressed through $E_{10}$, $E_{01}$, and $E_{11}$ by means of the following 
formulas:
$$
\pagebreak
\gather
\hskip -2em
E_{20}=\frac{1}{2}\,E_{10}^2-\frac{1}{2}\,L^2,
\kern 7em
E_{02}=\frac{1}{2}\,E_{01}^2-L^2,
\mytag{1.7}\\
\vspace{2ex}
\hskip -2em
\aligned
E_{21}&=\frac{2\,E_{10}^3\,E_{11}+2\,E_{01}^2\,E_{10}\,E_{11}
-E_{01}\,E_{10}^4+E_{01}^5}
{8\,(E_{01}^2+E_{10}^2)\vphantom{\vrule height 11pt}}\,+\\
\vspace{1ex}
&+\,\frac{6\,E_{10}\,E_{11}\,L^2-2\,E_{01}\,E_{10}^2\,L^2
-8\,E_{01}^3\,L^2+3\,E_{01}\,L^4}
{8\,(E_{01}^2+E_{10}^2)\vphantom{\vrule height 11pt}},
\endaligned
\mytag{1.8}\\
\vspace{2ex}
\hskip -2em
\aligned
E_{12}=&\frac{E_{01}^4\,E_{10}-2\,E_{01}^3\,E_{11}
-2\,E_{01}\,E_{10}^2\,E_{11}-E_{10}^5}
{8\,(E_{01}^2+E_{10}^2)\vphantom{\vrule height 11pt}}\,+\\
\vspace{1ex}
&\kern 5em +\,\frac{6\,E_{10}^3\,L^2-6\,E_{01}\,E_{11}\,L^2+3\,E_{10}\,L^4}
{8\,(E_{01}^2+E_{10}^2)\vphantom{\vrule height 11pt}},
\endaligned
\mytag{1.9}\\
\vspace{2ex}
\hskip -2em
E_{30}=-\frac{1}{3}\,E_{12}-\frac{1}{6}\,E_{10}\,E_{01}^2
-\frac{1}{2}\,E_{10}\,L^2+\frac{1}{6}\,E_{10}^3
+\frac{1}{3}\,E_{01}\,E_{11},\quad
\mytag{1.10}\\
\vspace{2ex}
\hskip -2em
E_{03}=-\frac{1}{3}\,E_{21}-\frac{1}{6}\,E_{01}\,E_{10}^2-\frac{5}{6}\,E_{01}\,L^2
+\frac{1}{6}\,E_{01}^3+\frac{1}{3}\,E_{10}\,E_{11}.\quad
\mytag{1.11}
\endgather
$$\par
     A general solution for the equation \mythetag{1.6} was derived in \mycite{49},
including a two-parameter solution and several one-parameter solutions. As it was
proved in \mycite{50}, the one-parameter solutions do not lead to perfect cuboids.
The two-parameter solution of \mythetag{1.6} is written in \mycite{49} for the case
of a rational cuboid with unit space diagonal $L=1$. This solution with the 
parameters $b$ and $c$ is given by the formulas 
$$
\align
&\hskip -2em
E_{11}=-\frac{b\,(c^2+2-4\,c)}{b^2\,c^2+2\,b^2-3\,b^2\,c+c-b\,c^2\,+2\,b},
\mytag{1.12}\\
\vspace{2ex}
&\hskip -2em
E_{01}=-\frac{b\,(c^2+2-2\,c)}{b^2\,c^2+2\,b^2-3\,b^2\,c+c-b\,c^2+2\,b},
\mytag{1.13}\\
\vspace{2ex}
&\hskip -2em
E_{10}=-\frac{b^2\,c^2+2\,b^2-3\,b^2\,c\,-c}{b^2\,c^2+2\,b^2-3\,b^2\,c
+c-b\,c^2+2\,b}.
\mytag{1.14}
\endalign
$$
Substituting \mythetag{1.12}, \mythetag{1.13}, and \mythetag{1.14} into the
formulas \mythetag{1.7}, \mythetag{1.8}, \mythetag{1.9}, 
\mythetag{1.10}, \mythetag{1.11} and taking into account that $L=1$, one can 
find that 
$$
\allowdisplaybreaks
\gather
\hskip -2em
\gathered
E_{20}=\frac{b}{2}\,(b\,c^2-2\,c-2\,b)\,(2\,b\,c^2-c^2-6\,b\,c+2
+4\,b)\,\times\\
\times\,(b\,c-1-b)^{-2}\,(b\,c-c-2\,b)^{-2},
\endgathered
\mytag{1.15}\\
\vspace{1ex}
\gathered
E_{02}=\frac{1}{2}\,(28\,b^2\,c^2-16\,b^2\,c-2\,c^2-4\,b^2-b^2\,c^4
+\,4\,b^3\,c^4-12\,b^3\,c^3\,+\\
+\,4\,b\,c^3+24\,b^3\,c-8\,b\,c-2\,b^4\,c^4
+12\,b^4\,c^3-26\,b^4\,c^2-8\,b^2\,c^3\,+\\
+\,24\,b^4\,c-16\,b^3-8\,b^4)\,(b\,c-1-b)^{-2}\,(b\,c-c-2\,b)^{-2}.
\endgathered\qquad\quad
\mytag{1.16}\\
\vspace{1ex}
\gathered
E_{21}=\frac{b}{2}\,(5\,c^6\,b-2\,c^6\,b^2+52\,c^5\,b^2-16\,c^5\,b
-2\,c^7\,b^2+2\,b^4\,c^8\,+\\
+\,142\,b^4\,c^6-26\,b^4\,c^7-426\,b^4\,c^5-61\,b^3\,c^6+100\,b^3\,c^5
+14\,c^7\,b^3\,-\\
-\,c^8\,b^3-20\,b\,c^2-8\,b^2\,c^2-16\,b^2\,c-128\,b^2\,c^4-200\,b^3\,c^3\,+\\
+\,244\,b^3\,c^2+32\,b\,c^3-112\,b^3\,c+768\,b^4\,c^4-852\,b^4\,c^3
+568\,b^4\,c^2\,+\\
+\,104\,b^2\,c^3-208\,b^4\,c+8\,c^4-4\,c^3+16\,b^3+32\,b^4-2\,c^5)\,\times\\
\times\,(b^2\,c^4-6\,b^2\,c^3+13\,b^2\,c^2-12\,b^2\,c+4\,b^2+c^2)^{-1}\,\times\\
\times\,(b\,c-1-b)^{-2}\,(b\,c-c-2\,b)^{-2},
\endgathered\qquad\quad
\mytag{1.17}\\
\vspace{2ex}
\hskip -2em
\gathered
E_{12}=(16\,b^6+32\,b^5-6\,c^5\,b^2+2\,c^5\,b-62\,b^5\,c^6+62\,b^6\,c^6\,-\\
-\,180\,b^6\,c^5+18\,b^5\,c^7-12\,b^6\,c^7-2\,b^5\,c^8+b^6\,c^8+248\,b^5\,c^2\,+\\
+\,248\,b^6\,c^2-96\,b^6\,c+321\,b^6\,c^4-180\,b^5\,c^3-144\,b^5\,c
-360\,b^6\,c^3\,+\\
+\,b^4\,c^8+8\,b^4\,c^6-6\,b^4\,c^7+18\,b^4\,c^5+7\,b^3\,c^6+90\,b^5\,c^5
-14\,b^3\,c^5\,-\\
-\,c^7\,b^3+17\,b^2\,c^4+28\,b^3\,c^3-28\,b^3\,c^2-4\,b\,c^3+8\,b^3\,c
-57\,b^4\,c^4\,+\\
+\,36\,b^4\,c^3+32\,b^4\,c^2-12\,b^2\,c^3-48\,b^4\,c-c^4+16\,b^4)\,\times\\
\times\,(b^2\,c^4-6\,b^2\,c^3+13\,b^2\,c^2-12\,b^2\,c+4\,b^2+c^2)^{-1}\,\times\\
\times\,(b\,c-1-b)^{-2}\,(b\,c-c-2\,b)^{-2}.
\endgathered\qquad
\mytag{1.18}\\
\gathered
E_{03}=\frac{b}{2}\,(b^2\,c^4-5\,b^2\,c^3+10\,b^2\,c^2-10\,b^2\,c+4\,b^2+2\,b\,c\,+\\
+\,2\,c^2-b\,c^3)\,(2\,b^2\,c^4-12\,b^2\,c^3+26\,b^2\,c^2-24\,b^2\,c\,+\\
+\,8\,b^2-c^4\,b+3\,b\,c^3-6\,b\,c+4\,b+c^3-2\,c^2+2\,c)\,\times\\
\times\,((b^2\,c^4-6\,b^2\,c^3+13\,b^2\,c^2-12\,b^2\,c+4\,b^2+c^2)^{-1}\,\times\\
\times\,(b\,c-1-b)^{-2}\,(-c+b\,c-2\,b)^{-2},
\endgathered\qquad\quad
\mytag{1.19}\\
\vspace{1ex}
\gathered
E_{30}=c\,b^2\,(1-c)\,(c-2)\,(b\,c^2-4\,b\,c+2+4\,b)\,\times\\
\times\,(2\,b\,c^2-c^2-4\,b\,c+2\,b)\,\times\\
\times\,(b^2\,c^4-6\,b^2\,c^3+13\,b^2\,c^2-12\,b^2\,c+4\,b^2+c^2)^{-1}\,\times\\
\times\,(b\,c-1-b)^{-2}\,(-c+b\,c-2\,b)^{-2}.
\endgathered\qquad\quad
\mytag{1.20}
\endgather
$$
The formulas \mythetag{1.15}, \mythetag{1.16}, \mythetag{1.17}, \mythetag{1.18}, 
\mythetag{1.19}, and \mythetag{1.20} were derived in \mycite{50}. Based on these 
formulas, two inverse problems were formulated. 
\myproblem{1.1} Find all pairs of rational numbers $b$ and $c$ for which the
cubic equations \mythetag{1.2} and \mythetag{1.3} with the coefficients 
\mythetag{1.14}, \mythetag{1.15},	\mythetag{1.20}, \mythetag{1.13}, 
\mythetag{1.16},	\mythetag{1.19} possess positive rational roots $x_1$, $x_2$, 
$x_3$, $d_1$, $d_2$, $d_3$ obeying the auxiliary equations \mythetag{1.5} with 
the right hand sides \mythetag{1.17}, \mythetag{1.18}, \mythetag{1.12}. 
\endproclaim
\myproblem{1.2} Find at least one pair of rational numbers $b$ and $c$ for which 
the cubic equations \mythetag{1.2} and \mythetag{1.3} with the coefficients 
\mythetag{1.14}, \mythetag{1.15},	\mythetag{1.20}, \mythetag{1.13}, 
\mythetag{1.16},	\mythetag{1.19} possess positive rational roots $x_1$, $x_2$, 
$x_3$, $d_1$, $d_2$, $d_3$ obeying the auxiliary equations \mythetag{1.5} with 
the right hand sides \mythetag{1.17}, \mythetag{1.18}, \mythetag{1.12}. 
\endproclaim
     The problems~\mytheproblem{1.1} and \mytheproblem{1.2} are equivalent to
finding all perfect cuboids and to finding at least one perfect cuboid respectively.
In the present paper we study the cubic equations \mythetag{1.2} and \mythetag{1.3}
for reducibility using the methods of \mycite{41}, which were applied to a twelfth
order Diophantine equation in that paper. 
\head
2. The first reducibility case $b=0$. 
\endhead
     Note that the formulas \mythetag{1.12} through \mythetag{1.20} for the 
coefficients of the cubic equations \mythetag{1.2} and \mythetag{1.3} and for the 
right hand sides of the auxiliary equations \mythetag{1.5} possess denominators. 
The simultaneous non-vanishing condition for all of their denominators is written 
as follows:
$$
\hskip -2em
\gathered
(b^2\,c^4-6\,b^2\,c^3+13\,b^2\,c^2-12\,b^2\,c+4\,b^2+c^2)\,\times\\
\times\,(b\,c-1-b)\,(b\,c-c-2\,b)\neq 0.
\endgathered
\mytag{2.1}
$$
The case $b=0$ is very simple. The non-vanishing condition \mythetag{2.1} in this 
case is written as $c\neq 0$. Substituting $b=0$ into \mythetag{1.14}, \mythetag{1.15},	
\mythetag{1.20}, \mythetag{1.13}, \mythetag{1.16},	\mythetag{1.19}, we find that 
the cubic equations \mythetag{1.2} and \mythetag{1.3} reduce to 
$$
\xalignat 2
&\hskip -2em
x^2\,(x-1)=0,
&&d\,(d-1)\,(d+1)=0.
\mytag{2.2}
\endxalignat
$$
Substituting $b=0$ into \mythetag{1.17}, \mythetag{1.18}, \mythetag{1.12}, we obtain
$$
\xalignat 3
&\hskip -2em
E_{21}=0,
&&E_{11}=0,
&&E_{12}=-1.
\mytag{2.3}
\endxalignat
$$
\par
     The equations \mythetag{2.2} are already factored. We can choose their roots 
as follows:
$$
\xalignat 3
&\hskip -2em
x_1=1,
&&x_2=0,
&&x_3=0,\\
\vspace{-1.7ex}
\mytag{2.4}\\
\vspace{-1.7ex}
&\hskip -2em
d_1=0,
&&d_1=1,
&&d_3=-1.
\endxalignat
$$
Substituting \mythetag{2.3} and \mythetag{2.4} into \mythetag{1.5}, we find that
the auxiliary equations \mythetag{1.5} are fulfilled. Substituting \mythetag{2.4}
along with $L=1$ into \mythetag{1.1}, we find that the cuboid 
equations \mythetag{1.1} are also fulfilled. However, the formulas \mythetag{2.4} 
do not provide a perfect cuboid since its edges and face diagonals cannot be zero 
or negative.\par
\mytheorem{2.1} If\/ $b=0$ and $c\neq 0$, then the cubic polynomials in \mythetag{1.2}
and \mythetag{1.3} are reducible and provide three integer roots for each of the 
equations \mythetag{1.2} and \mythetag{1.3} satisfying the auxiliary equations 
\mythetag{1.5} but not resolving the problem~\mytheproblem{1.2}.
\endproclaim
\head
3. The second reducibility case $c=0$. 
\endhead
     The case $c=0$ is also simple. The non-vanishing condition \mythetag{2.1} in this
case turns to $b\,(1+b)\neq 0$. The cubic equations \mythetag{1.2} and \mythetag{1.3}
turn to 
$$
\hskip -2em
\aligned
x\,&\bigl(2\,(1+b)^2\,x^2+2\,b\,(1+b)\,x-(1+2\,b)\bigr)=0,\\
\vspace{1ex}
(d+1)\,&\bigl(2\,(1+b)^2\,d^{\kern 1pt 2}-2\,b\,(1+b)\,d-(1+2\,b)\bigr)=0.
\endaligned
\mytag{3.1}
$$
Substituting $c=0$ into the formulas \mythetag{1.17}, \mythetag{1.18}, and 
\mythetag{1.12}, we obtain
$$
\xalignat 3
&\hskip -2em
E_{21}=\frac{1+2\,b}{2\,(1+b)^2},
&&E_{11}=\frac{-1}{1+b},
&&E_{12}=1.
\mytag{3.2}
\endxalignat
$$
\par
     The equations \mythetag{3.1} are already factored. Upon splitting off the linear 
terms they turn to quadratic equations very similar to each other:
$$
\hskip -2em
\aligned
&2\,(1+b)^2\,x^2+2\,b\,(1+b)\,x-(1+2\,b)=0,\\
\vspace{1ex}
&2\,(1+b)^2\,d^{\kern 1pt 2}-2\,b\,(1+b)\,d-(1+2\,b)=0.
\endaligned
\mytag{3.3}
$$
The discriminants of the quadratic equations \mythetag{3.3} do coincide:
$$
\hskip -2em
D=4\,(b^2+4\,b+2)\,(1+b)^2.
\mytag{3.4}
$$
The formula \mythetag{3.4} means that in order to solve the equations \mythetag{3.3}
in rational numbers one should resolve the following quadratic equation in rational 
numbers:
$$
\hskip -2em
b^2+4\,b+2=\beta^2.
\mytag{3.5}
$$
The equation \mythetag{3.5} can be written as $(b+2)^2-\beta^2=2$. The lemma~2.2 from
\mycite{49} is applicable to this equation. Applying this lemma, we get
$$
\xalignat 2
&\hskip -2em
b=\frac{t^2-4\,t+2}{2\,t},
&&\beta=\frac{t^2-2}{2\,t}
\mytag{3.6}
\endxalignat
$$
for some rational $t\neq 0$. Note that the condition $b\,(1+b)\neq 0$ derived from 
\mythetag{2.1} for $c=0$ implies no restrictions for $t$ since the equations
$$
\xalignat 2
&\frac{t^2-4\,t+2}{2\,t}=0,
&&\frac{t^2-4\,t+2}{2\,t}=-1
\endxalignat
$$
have no rational roots. Now, substituting $b$ from \mythetag{3.6} into the
quadratic equations \mythetag{3.3}, we can find their roots. Since the roots 
$x=0$ and $d=-1$ of the cubic equations \mythetag{3.1} are already known, we 
can write formulas for all their roots:
$$
\xalignat 3
&\hskip -2em
x_1=0,
&&x_2=\frac{-t\,(t-2)}{(t-1)^2+1},
&&x_3=\frac{2\,(t-1)}{(t-1)^2+1},
\qquad\\
\vspace{-1.3ex}
\mytag{3.7}\\
\vspace{-1.3ex}
&\hskip -2em
d_1=-1,
&&d_2=\frac{-2\,(t-1)}{(t-1)^2+1},
&&d_3=\frac{t\,(t-2)}{(t-1)^2+1}.
\qquad
\endxalignat
$$
Taking into account the formula for $b$ in \mythetag{3.6}, we transform 
\mythetag{3.2} to
$$
\xalignat 3
&\hskip -2em
E_{21}=\frac{2\,t\,(t^2-3\,t+2)}{(t-1)^2+1},
&&E_{11}=\frac{-2\,t}{(t-1)^2+1},
&&E_{12}=1.\quad
\mytag{3.8}
\endxalignat
$$
Now, if we substitute the formulas \mythetag{3.7} and \mythetag{3.8} into
\mythetag{1.5}, we find that the auxiliary equations \mythetag{1.5} are fulfilled. 
Similarly, substituting \mythetag{3.7} along with $L=1$ into \mythetag{1.1}, we 
find that the cuboid equations \mythetag{1.1} are also fulfilled. However, the 
formulas \mythetag{3.7} do not provide a perfect cuboid since $x_1$ is zero and
$d_1$ is negative.
\mytheorem{3.1} If\/ $c=0$ and $b\,(1+b)\neq 0$, then the cubic polynomials in 
\mythetag{1.2} and \mythetag{1.3} are reducible over $\Bbb Q$. Moreover, if\/ $b$ 
is given by the first formula \mythetag{3.6} for some rational $t\neq 0$, then
each of the cubic equations \mythetag{1.2} and \mythetag{1.3} has three rational 
roots satisfying the auxiliary equations \mythetag{1.5} but not resolving the 
problem~\mytheproblem{1.2}.
\endproclaim
\head
4. The third reducibility case $c=1$. 
\endhead
     The case $c=1$ is similar to the previous case $c=0$. The non-vanishing 
condition \mythetag{2.1} in this case turns to $b\neq -1$. The cubic 
equations \mythetag{1.2} and \mythetag{1.3} turn to 
$$
\hskip -2em
\aligned
x\,&\bigl(2\,(1+b)^2\,x^2-(2\,b+2)\,x-b\,(b+2)\bigr)=0,\\
\vspace{1ex}
(d+1)\,&\bigl(2\,(1+b)^2\,d^{\kern 1pt 2}-(2\,b+2)\,d-b\,(b+2)\bigr)=0.
\endaligned
\mytag{4.1}
$$
Substituting $c=1$ into the formulas \mythetag{1.17}, \mythetag{1.18}, and 
\mythetag{1.12}, we obtain
$$
\xalignat 3
&\hskip -2em
E_{21}=\frac{b\,(b+2)}{2\,(1+b)^2},
&&E_{11}=\frac{b}{1+b},
&&E_{12}=-1.
\mytag{4.2}
\endxalignat
$$
\par
     The equations \mythetag{4.1} are already factored. Upon splitting out linear 
terms they turn to quadratic equations. These two quadratic equations do coincide 
with each other up to the change of $x$ for $d$ and vice versa:
$$
\hskip -2em
\aligned
&2\,(1+b)^2\,x^2-(2\,b+2)\,x-b\,(b+2)=0,\\
\vspace{1ex}
&2\,(1+b)^2\,d^{\kern 1pt 2}-(2\,b+2)\,d-b\,(b+2)=0.
\endaligned
\mytag{4.3}
$$
Certainly, the discriminants of the coinciding equations \mythetag{4.3} do also 
coincide:
$$
\hskip -2em
D=4\,(2\,b^2+4\,b+1)\,(b+1)^2.
\mytag{4.4}
$$
The formula \mythetag{4.4} means that in order to solve the equations \mythetag{4.3}
in rational numbers one should resolve the following quadratic equation in rational 
numbers:
$$
\hskip -2em
2\,b^2+4\,b+1=\beta^2.
\mytag{4.5}
$$
If $b=0$, then \mythetag{4.5} implies $\beta=\pm\,1$, which yields two trivial 
rational solutions for the equation \mythetag{4.5}. If $b\neq 0$, we can write
the equation \mythetag{4.5} as 
$$
\hskip -2em
\frac{1}{b^{\kern 1pt 2}}+\frac{4}{b}+2=\frac{\beta^{\kern 1pt 2}}{b^{\kern 1pt 2}}.
\mytag{4.6}
$$
The equation \mythetag{4.6} is very similar to \mythetag{3.5}. Therefore the
formulas \mythetag{3.6} yield 
$$
\xalignat 2
&\hskip -2em
\frac{1}{b}=\frac{t^2-4\,t+2}{2\,t},
&&\frac{\beta}{b}=\frac{t^2-2}{2\,t}
\mytag{4.7}
\endxalignat
$$
for some rational $t\neq 0$. Then the formulas \mythetag{4.7} can be transformed to
$$
\xalignat 2
&\hskip -2em
b=\frac{2\,t}{(t-2)^2-2},
&&\beta=\frac{t^2-2}{(t-2)^2-2}. 
\mytag{4.8}
\endxalignat
$$
Note that the formulas \mythetag{4.8} are consistent since their denominators 
cannot vanish for any rational $t$. Note also that one of the above trivial 
solutions with $b=0$ and $\beta=-1$ can be obtained from \mythetag{4.8} for $t=0$. 
In order to cover the other trivial solution with $b=0$ and $\beta=1$ we need to 
add the following sign option to \mythetag{4.8}:
$$
\xalignat 2
&\hskip -2em
b=\frac{2\,t}{(t-2)^2-2},
&&\beta=\pm\,\frac{t^2-2}{(t-2)^2-2}. 
\mytag{4.9}
\endxalignat
$$
Due to these observations the restriction $t\neq 0$ is removed and the formulas 
\mythetag{4.9} cover all rational solutions of the equation \mythetag{4.5}. 
The condition $b\neq -1$ derived from \mythetag{2.1} for $c=1$ implies no
restrictions for $t$ since the equation
$$
\frac{2\,t}{(t-2)^2-2}=-1
$$
has no rational roots. Now we can substitute $b$ from \mythetag{4.9} into the 
quadratic equations \mythetag{4.3} and find their roots. Since the roots $x=0$ 
and $d=-1$ of the cubic equations \mythetag{4.1} are already known, we can write 
formulas for all their roots:
$$
\pagebreak
\xalignat 3
&\hskip -2em
x_1=0,
&&x_2=\frac{t\,(t-2)}{(t-1)^2+1},
&&x_3=\frac{-2\,(t-1)}{(t-1)^2+1},
\qquad\\
\vspace{-1.3ex}
\mytag{4.10}\\
\vspace{-1.3ex}
&\hskip -2em
d_1=-1,
&&d_2=\frac{-2\,(t-1)}{(t-1)^2+1},
&&d_3=\frac{t\,(t-2)}{(t-1)^2+1}.
\qquad
\endxalignat
$$
Taking into account the formula for $b$ in \mythetag{4.9}, we transform 
\mythetag{4.2} to
$$
\xalignat 3
&\hskip -2em
E_{21}=\frac{2\,(t^2-3\,t+2)\,t}{((t-1)^2+1)^2},
&&E_{11}=\frac{2\,t}{(t-1)^2+1},
&&E_{12}=-1.\qquad
\mytag{4.11}
\endxalignat
$$
Now, if we substitute the formulas \mythetag{4.10} and \mythetag{4.11} into
\mythetag{1.5}, we find that the auxiliary equations \mythetag{1.5} are fulfilled. 
Similarly, substituting \mythetag{4.10} along with $L=1$ into \mythetag{1.1}, we 
find that the cuboid equations are also fulfilled. But again, the formulas 
\mythetag{4.10} do not provide a perfect cuboid since $x_1$ is zero and $d_1$ is 
negative.\par
\mytheorem{4.1} If\/ $c=1$ and $b\neq -1$, then the cubic polynomials in 
\mythetag{1.2} and \mythetag{1.3} are reducible over $\Bbb Q$. Moreover, if\/ $b$ 
is given by the first formula \mythetag{4.9} for some rational $t$, then
each of the cubic equations \mythetag{1.2} and \mythetag{1.3} has three rational 
roots satisfying the auxiliary equations \mythetag{1.5} but not resolving the 
problem~\mytheproblem{1.2}.
\endproclaim
\head
5. The fourth reducibility case $c=2$. 
\endhead
     The case $c=2$ is similar to both of the previous cases $c=0$ and $c=1$. 
The non-vanishing condition \mythetag{2.1} in this case turns to $b\neq 1$. 
The cubic equations \mythetag{1.2} and \mythetag{1.3} in this case turn to 
the following ones:
$$
\hskip -2em
\aligned
x\,&\bigl(2\,(b-1)^2\,x^2+2\,(b-1)\,x-b\,(b-2)\bigr)=0,\\
\vspace{1ex}
(d+1)\,&\bigl(2\,(b-1)^2\,d^{\kern 1pt 2}-2\,(b-1)\,d-b\,(b-2)\bigr)=0.
\endaligned
\mytag{5.1}
$$
Substituting $c=2$ into the formulas \mythetag{1.17}, \mythetag{1.18}, and 
\mythetag{1.12}, we obtain
$$
\xalignat 3
&\hskip -2em
E_{21}=-\frac{b\,(b-2)}{2\,(b-1)^2},
&&E_{11}=-\frac{b}{b-1},
&&E_{12}=-1.
\mytag{5.2}
\endxalignat
$$\par
     The equations \mythetag{5.1} are already factored. If we split off the linear 
terms, they turn to quadratic equations. These two quadratic equations are
very similar:
$$
\hskip -2em
\aligned
&2\,(b-1)^2\,x^2+2\,(b-1)\,x-b\,(b-2)=0,\\
\vspace{1ex}
&2\,(b-1)^2\,d^{\kern 1pt 2}-2\,(b-1)\,d-b\,(b-2)=0.
\endaligned
\mytag{5.3}
$$
The discriminants of the quadratic equations \mythetag{5.3} do coincide:
$$
\hskip -2em
D=4\,(2\,b^2-4\,b+1)\,(b-1)^2.
\mytag{5.4}
$$
The formula \mythetag{5.4} means that in order to solve the equations \mythetag{5.3}
in rational numbers one should resolve the following quadratic equation in rational 
numbers:
$$
\hskip -2em
2\,b^2-4\,b+1=\beta^2.
\mytag{5.5}
$$
The equation \mythetag{5.5} is similar to the equation \mythetag{4.5}. If $b=0$, 
then it has two trivial solutions with $\beta=\pm\,1$. If $b\neq 0$, we can write
the equation \mythetag{5.5} as
$$
\hskip -2em
\biggl(\frac{1}{b}-2\biggr)^2\!\!-\frac{\beta^2}{b^2}=2.
\mytag{5.6}
$$
The lemma~2.2 from \mycite{49} is applicable to the equation \mythetag{5.6}.
This lemma yields
$$
\xalignat 2
&\hskip -2em
\frac{1}{b}=\frac{t^2+4\,t+2}{2\,t},
&&\frac{\beta}{b}=\frac{t^2-2}{2\,t}
\mytag{5.7}
\endxalignat
$$
for some rational $t\neq 0$. Now the formulas \mythetag{5.7} can be transformed to
$$
\xalignat 2
&\hskip -2em
b=\frac{2\,t}{(t+2)^2-2},
&&\beta=\frac{t^2-2}{(t+2)^2-2}. 
\mytag{5.8}
\endxalignat
$$
Like \mythetag{4.8}, the formulas \mythetag{5.8} are consistent since their 
denominators cannot vanish for any rational $t$. One of the above trivial solutions 
with $b=0$ and $\beta=-1$ can be obtained from \mythetag{5.8} for $t=0$. In order 
to cover the other trivial solution with $b=0$ and $\beta=1$ we need to add the
following sign option to \mythetag{5.8}:
$$
\xalignat 2
&\hskip -2em
b=\frac{2\,t}{(t+2)^2-2},
&&\beta=\pm\,\frac{t^2-2}{(t+2)^2-2}. 
\mytag{5.9}
\endxalignat
$$
Due to the above observations the restriction $t\neq 0$ is removed and the formulas 
\mythetag{5.9} cover all rational solutions of the equation \mythetag{5.5}. The 
condition $b-1\neq 0$ derived from \mythetag{2.1} for $c=2$ implies no restrictions 
for $t$ since the equation
$$
\frac{2\,t}{(t+2)^2-2}=1
$$
has no rational roots. Now we can substitute $b$ from \mythetag{5.9} into the 
quadratic equations \mythetag{5.3} and find their roots. Since the roots $x=0$ 
and $d=1$ of the cubic equations \mythetag{5.1} are already known, we can write 
formulas for all their roots:
$$
\xalignat 3
&\hskip -2em
x_1=0,
&&x_2=\frac{2\,(t+1)}{(t+1)^2+1},
&&x_3=\frac{t\,(t+2)}{(t+1)^2+1},
\qquad\\
\vspace{-1.3ex}
\mytag{5.10}\\
\vspace{-1.3ex}
&\hskip -2em
d_1=1,
&&d_2=\frac{-t\,(t+2)}{(t+1)^2+1},
&&d_3=\frac{-2\,(t+1)}{(t+1)^2+1}.
\qquad
\endxalignat
$$
Taking into account the formula for $b$ in \mythetag{5.9}, the formula 
\mythetag{5.2} is transformed to
$$
\xalignat 3
&\hskip -2em
E_{21}=\frac{2\,(t^2+3\,t+2)\,t}{((t+1)^2+1)^2},
&&E_{11}=\frac{2\,t}{(t+1)^2+1},
&&E_{12}=-1.\qquad
\mytag{5.11}
\endxalignat
$$
If we substitute the formulas \mythetag{5.10} and \mythetag{5.11} into
\mythetag{1.5}, we find that the auxiliary equations \mythetag{1.5} are fulfilled. 
Similarly, substituting \mythetag{5.10} into \mythetag{1.1} along with $L=1$, we 
find that the cuboid equations are also fulfilled. But, like in the previous cases, 
the formulas \mythetag{5.10} do not provide a perfect cuboid since $x_1$ is zero.
\par
\mytheorem{5.1} If\/ $c=2$ and $b\neq 1$, then the cubic polynomials in 
\mythetag{1.2} and \mythetag{1.3} are reducible over $\Bbb Q$. Moreover, if\/ $b$ 
is given by the first formula \mythetag{5.9} for some rational $t$, then
each of the cubic equations \mythetag{1.2} and \mythetag{1.3} has three rational 
roots satisfying the auxiliary equations \mythetag{1.5} but not resolving the 
problem~\mytheproblem{1.2}.
\endproclaim
\head
6. Some other reducibility cases. 
\endhead
     Note that in each of the previous four cases the cubic equation \mythetag{1.2}
has the root $x=0$. The necessary and sufficient condition for that is written as
$$
\hskip -2em
P(x)=E_{30}=0,
\mytag{6.1}
$$
where $P(x)$ is the cubic polynomial in the left hand side of the equation 
\mythetag{1.2}. Looking at the formulas \mythetag{6.1} and \mythetag{1.20}, we 
see that along with the conditions $b=0$, $c=0$, $c=1$, $c=2$, which were already 
considered in the previous cases, there are the following two conditions for 
vanishing $E_{30}$:
$$
\align
&b\,c^2-4\,b\,c+4\,b+2=0,
\mytag{6.2}\\
&2\,b\,c^2-4\,b\,c+2\,b-c^2=0.
\mytag{6.3}
\endalign
$$
Let's denoter through $Q(d)$ the cubic polynomial in the left hand side of the
equation \mythetag{1.3}. Then one can easily derive the following formulas:
$$
\gather
\hskip -2em
Q(-1)=-(c-1)^2\,(b\,c^2-4\,b\,c+4\,b+2)^2\,b^2\,c^2\,\times\\
\hskip -2em
\times\,(b^2\,c^4-6\,b^2\,c^3+13\,b^2\,c^2-12\,b^2\,c+4\,b^2+c^2)^{-1}
\,\times
\mytag{6.4}\\
\hskip -2em
\times\,(b\,c-1-b)^{-2}\,(b\,c-c-2\,b)^{-2},\\
\vspace{1ex}
\hskip -2em
Q(1)=(c-2)^2\,(2\,b\,c^2-4\,b\,c+2\,b-c^2)^2\,b^2\,\times\\
\hskip -2em
\times\,(b^2\,c^4-6\,b^2\,c^3+13\,b^2\,c^2-12\,b^2\,c+4\,b^2+c^2)^{-1}
\,\times
\mytag{6.5}\\
\hskip -2em
\times\,(b\,c-1-b)^{-2}\,(b\,c-c-2\,b)^{-2}.
\endgather
$$
Comparing \mythetag{6.2} with \mythetag{6.4}, we see that the condition
\mythetag{6.2} implies $Q(-1)=0$. Similarly, comparing \mythetag{6.3} with 
\mythetag{6.5}, we find that the condition \mythetag{6.3} implies $Q(1)=0$.
These observations yield the following theorems.
\mytheorem{6.1} If the condition \mythetag{6.2} is fulfilled, then both
cubic equations \mythetag{1.2} and \mythetag{1.3} are reducible. In this case 
the first of them has the root $x=0$, while the second has the root $d=-1$.
\endproclaim
\mytheorem{6.2} If the condition \mythetag{6.3} is fulfilled, then both
cubic equations \mythetag{1.2} and \mythetag{1.3} are reducible. In this case 
the first of them has the root $x=0$, while the second has the root $d=1$.
\endproclaim
\head
7. The fifth reducibility case. 
\endhead
    The fifth reducibility case is defined by the condition \mythetag{6.2}.
The equality \mythetag{6.2} is linear with respect to $b$. It can be resolved
as follows:
$$
\hskip -2em
b=\frac{-2}{(c-2)^2}.
\mytag{7.1}
$$
Substituting \mythetag{7.1} into the formulas \mythetag{1.14}, \mythetag{1.15}, 
and \mythetag{1.20}, we find that the first cubic equation \mythetag{1.2} reduces 
to the following one: 
$$
\hskip -2em
x\,\bigl((c^2-2\,c+2)\,x+2\,(c-1)\bigr)\bigl((c^2-2\,c+2)\,x-c\,(c-2)\bigr)=0. 
\mytag{7.2}
$$
It is easy to see that the equation \mythetag{7.2} has three rational roots
$$
\xalignat 3
&\hskip -2em
x_1=0,
&&x_2=\frac{-2\,(c-1)}{(c-1)^2+1},
&&x_3=\frac{c\,(c-2)}{(c-1)^2+1}.
\qquad
\mytag{7.3}
\endxalignat
$$
Now let's substitute \mythetag{7.1} into the formulas \mythetag{1.13}, \mythetag{1.16}, 
and \mythetag{1.19}. As a result we find that the second cubic equation \mythetag{1.3} 
reduces to the following one: 
$$
\hskip -2em
(d+1)\,\bigl((c^2-2\,c+2)\,d-2\,(c-1)\bigr)
\bigl((c^2-2\,c+2)\,d-c\,(c-2)  \bigr)=0. 
\mytag{7.4}
$$
Again, it is easy to see that the equation \mythetag{7.4} has three rational roots
$$
\xalignat 3
&\hskip -2em
d_1=-1,
&&d_2=\frac{c\,(c-2)}{(c-1)^2+1},
&&d_3=\frac{2\,(c-1)}{(c-1)^2+1}.
\qquad
\mytag{7.5}
\endxalignat
$$
Now let's substitute \mythetag{7.1} into the formulas \mythetag{1.17}, 
\mythetag{1.18}, and \mythetag{1.12}. This yields
$$
\align
&\hskip -2em
E_{21}=\frac{2\,(c-2)\,(c-1)\,c}{((c-1)^2+1)^2},\\
\vspace{1ex}
&\hskip -2em
E_{11}=\frac{2\,((c-2)^2-2)\,(c-2)}{((c-1)^2+1)^2},
\mytag{7.6}\\
\vspace{1ex}
&\hskip -2em
E_{12}=\frac{-((c-2)^2-2)\,(c^2-2)}{((c-1)^2+1)^2}.
\endalign
$$
If we substitute the formulas \mythetag{7.3}, \mythetag{7.5}, and \mythetag{7.6}
into \mythetag{1.5}, we find that the auxiliary equations \mythetag{1.5} are 
fulfilled. Similarly, substituting \mythetag{7.3} and \mythetag{7.5} along with 
$L=1$ into the equations \mythetag{1.1}, we find that the cuboid equations are 
also fulfilled. But, like in the previous cases, the formulas \mythetag{7.3} and 
\mythetag{7.5} do not provide a perfect cuboid since $x_1$ is zero.\par
     Note that \mythetag{7.1} provides the restriction $c\neq 2$. Substituting
\mythetag{7.1} into \mythetag{2.1}, we get the condition $(c-1)^2+1\neq 0$ which
is fulfilled for all rational $c$. Summarizing the results of this section, we
can formulate the following theorem.
\mytheorem{7.1} If\/ $c\neq 2$ and\/ $b\,(c-2)^2=-2$, then the cubic polynomials 
in \mythetag{1.2} and \mythetag{1.3} are reducible over the field of rational 
numbers $\Bbb Q$. Moreover, each of the cubic equations \mythetag{1.2} and 
\mythetag{1.3} has three rational roots satisfying the auxiliary equations 
\mythetag{1.5} but not resolving the problem~\mytheproblem{1.2}.
\endproclaim
\head
8. The sixth reducibility case. 
\endhead
    The sixth reducibility case is defined by the condition \mythetag{6.3}.
The equality \mythetag{6.3} is linear with respect to $b$. It can be resolved
as follows:
$$
\hskip -2em
b=\frac{ c^2}{2\,(c-1)^2}.
\mytag{8.1}
$$
Substituting \mythetag{8.1} into the formulas \mythetag{1.14}, \mythetag{1.15}, 
and \mythetag{1.20}, we find that the first cubic equation \mythetag{1.2} reduces 
to the following one: 
$$
\hskip -2em
x\,\bigl((c^2-2\,c+2)\,x+2\,(c-1)\bigr)\bigl((c^2-2\,c+2)\,x-c\,(c-2)\bigr)=0. 
\mytag{8.2}
$$
It is easy to see that the equation \mythetag{8.2} has three rational roots
$$
\xalignat 3
&\hskip -2em
x_1=0,
&&x_2=\frac{-2\,(c-1)}{(c-1)^2+1)}
&&x_3=\frac{c\,(c-2)}{(c-1)^2+1}.
\qquad
\mytag{8.3}
\endxalignat
$$
Now let's substitute \mythetag{8.1} into the formulas \mythetag{1.13}, \mythetag{1.16}, 
and \mythetag{1.19}. As a result we find that the second cubic equation \mythetag{1.3} 
reduces to the following one: 
$$
\hskip -2em
(d-1)\,\bigl((c^2-2\,c+2)\,d-2\,(c-1)\bigr)
\bigl((c^2-2\,c+2)\,d-c\,(c-2)  \bigr)=0. 
\mytag{8.4}
$$
Again, it is easy to see that the equation \mythetag{8.4} has three rational roots
$$
\xalignat 3
&\hskip -2em
d_1=1,
&&d_2=\frac{c\,(c-2)}{(c-1)^2+1},
&&d_3=\frac{2\,(c-1)}{(c-1)^2+1}.
\qquad
\mytag{8.5}
\endxalignat
$$
The formulas \mythetag{8.3} and \mythetag{8.5} are almost identical to the formulas
\mythetag{7.3} and \mythetag{7.5}. The only difference is the sign of $d_1$.\par
     Now let's substitute \mythetag{7.1} into the formulas \mythetag{1.17}, 
\mythetag{1.18}, and \mythetag{1.12}. This yields
$$
\align
&\hskip -2em
E_{21}=-\frac{2\,(c-2)\,(c-1)\,c}{((c-1)^2+1)^2},\\
\vspace{1ex}
&\hskip -2em
E_{11}=\frac{2\,c\,((c-2)^2-2)\,(c-1)}{((c-1)^2+1)^2},
\mytag{8.6}\\
\vspace{1ex}
&\hskip -2em
E_{12}=\frac{((c-2)^2-2)\,(c^2-2)}{((c-1)^2+1)^2}.
\endalign
$$
If we substitute the formulas \mythetag{8.3}, \mythetag{8.5}, and \mythetag{8.6}
into \mythetag{1.5}, we find that the auxiliary equations \mythetag{1.5} are 
fulfilled. Similarly, substituting \mythetag{8.3} and \mythetag{8.5} along with 
$L=1$ into the equations \mythetag{1.1}, we find that the cuboid equations are 
also fulfilled. But like in all previous cases, the formulas 
\mythetag{8.3} and \mythetag{8.5} do not provide a perfect cuboid since $x_1$ 
is zero.\par
     Note that \mythetag{8.1} provides the restriction $c\neq 1$. Substituting
\mythetag{7.1} into \mythetag{2.1}, we get the condition $c\,((c-2)^2+1)\neq 0$,
which reduces to $c\neq 0$ since $(c-2)^2+1$ is always positive. Summarizing the 
results of this section, we can formulate the following theorem.
\mytheorem{8.1} If\/ $c\neq 0$, $c\neq 1$, and\/ $2\,b\,(c-1)^2=c^2$, then the 
cubic polynomials in \mythetag{1.2} and \mythetag{1.3} are reducible over the 
field of rational numbers $\Bbb Q$. Moreover, each of the cubic equations 
\mythetag{1.2} and \mythetag{1.3} has three rational roots satisfying the auxiliary 
equations \mythetag{1.5} but not resolving the problem~\mytheproblem{1.2}.
\endproclaim
\head
9. Concluding remarks.
\endhead
     Six special cases of reducibility considered in the above sections do not
exhaust all of the cases where the cubic equations \mythetag{1.2} and \mythetag{1.3}
are reducible over $\Bbb Q$. There is one very special case with the following values
of $b$ and $c$:
$$
\pagebreak
\xalignat 2
&b=\frac{14}{5},
&&c=-\frac{7}{2}.
\endxalignat
$$
In this very special case the equations \mythetag{1.2} and \mythetag{1.3} are factored 
as follows:
$$
\hskip -2em
\aligned
&(17\ x+15)\,(9248\ x^2+3128\ x-495)=0,\\
&(17\ d+8)\,(9248\ d^{\kern 1pt 2}-952\ d-8175)=0.
\endaligned
\mytag{9.1}
$$
The search for such special cases of reducibility like \mythetag{9.1} is planned for 
the future.\par

\enddocument
\end